\documentclass[reqno,11pt]{amsart}

\usepackage[all]{xy}

\numberwithin{equation}{section}

\theoremstyle{plain}
\newtheorem{lemma}{Lemma}[section]
\newtheorem{theor}[lemma]{Theorem}
\newtheorem{prop}[lemma]{Proposition}

\newtheorem{claim}[lemma]{Claim}

\theoremstyle{definition}

\theoremstyle{remark}
\newtheorem{rem}[lemma]{Remark}
\newtheorem{nots}[lemma]{Notation}

\newcommand{\pic}{{\rm Pic\thinspace}}

\newcommand{\vdim}{{\rm vdim\thinspace}}
\newcommand{\edim}{{\rm edim\thinspace}}

\newcommand{\p}{\mathbb{P}}

\newcommand{\z}{\mathbb{Z}}

\newcommand{\n}{\mathbb{N}}

\newcommand{\lh}{\hat{l}}
\newcommand{\vh}{\hat{v}}

\newcommand{\co}{{\mathcal O}}
\newcommand{\cl}{{\mathcal L}}
\newcommand{\clh}{\widehat{{\mathcal L}}}
\newcommand{\cm}{{\mathcal M}}

\newcommand{\ci}{{\mathcal I}}
\newcommand{\crsi}{{\mathcal R}_{\st}}
\newcommand{\crpi}{{\mathcal R}_{\cup L_i}}
\newcommand{\ccr}{{\mathcal R}}
\newcommand{\cx}{{\mathcal X}}

\newcommand{\fb}{\mathfrak{b}}

\newcommand{\st}{\widetilde{S}}
\newcommand{\ct}{\widetilde{C}}

\newcommand{\pp}{{\mathbf P}}

\newcommand{\svs}{\vspace{1mm}}

\newcommand{\ind}{\rule{2mm}{0mm}}

\newcommand{\map}{\rightarrow}

\newcommand{\vrk}{\hfill $\Box$}

\begin{document}

\title[Degeneration of linear systems on $K3$ surfaces]{Degeneration of linear systems through fat points on $K3$ surfaces}
\author{Cindy De Volder}
\address{
Department of Pure Mathematics and Computeralgebra,
Galglaan 2, \newline B-9000 Ghent, Belgium}
\email{cdv@cage.ugent.be}
\thanks{The first author is a Postdoctoral Fellow of
the Fund for Scientific Research-Flanders (Belgium) (F.W.O.-Vlaanderen)}

\author{Antonio Laface}
\address{
Dipartimento di Matematica, Universit\`a degli Studi di Milano,
Via Saldini 50, \newline 20100 Milano, Italy
}
\email{antonio.laface@unimi.it}
\thanks{The second author would like to thank the European Research and
Training Network EAGER for the support provided at Ghent
University. He also acknowledges the support of the MIUR of the
Italian Government in the framework of the National Research
Project ``Geometry in Algebraic Varieties'' (Cofin 2002)}
\keywords{Linear systems, fat points, generic $K3$ surfaces.}
\subjclass{14C20, 14J28.}
\begin{abstract}
In this paper we introduce a technique to degenerate $K3$ surfaces
and linear systems through fat points in general position on $K3$ surfaces.
Using this degeneration we show that on generic $K3$ surfaces
it is enough to prove that linear systems with one fat point are non-special
in order to obtain the non-speciality of homogeneous linear systems through
$n = 4^u9^w$ fat points in general position.
Moreover, we use this degeneration to obtain a result
for homogeneous linear systems through
$n = 4^u9^w$ fat points in general position
on a general quartic surface in $\p^3$.
\end{abstract}
\maketitle

\section{Introduction}

In this paper we assume the ground field is the field of the complex numbers.

Let $S$ be a smooth projective generic $K3$ surface (i.e. $\pic(S) \cong \z$)
and let $H$ be the generator of $\pic(S)$.

Consider $n$ points in general position on $S$, to each one of them associate a natural number $m_i$
called the {\em multiplicity} of the point and
let $n_j$ be the number of $p_i$ with multiplicity $m_i$.

For a linear system
of curves in $|dH|$ with $n_j$ general base points of
multiplicity $m_j$ for $j=1\cdots k$,
define its {\em virtual dimension} $v$ as $\dim|dH|-\sum n_i m_i(m_i+1)/2$ and its
{\em expected dimension} by $e=\max\{v,-1\}$.
If the dimension of the linear system is $l$,
then
$ v \leq e \leq l$.

Observe that it is possible to have $e < l$,
since the
conditions imposed by the points may be dependent.
In this case we say that the system is {\em special}.

Linear systems through general fat points on rational surfaces
have been studied by many authors
(see e.g. \cite{AG,BH,CM1,CM2}),
but, as far as we know, on $K3$ surfaces,
no results on the non-speciality of such systems are known.

In section~\ref{degen} we develop a technique to degenerate a $K3$ surface
and linear systems through fat points on $K3$ surfaces,
this degeneration is similar to the degeneration of the
plane used in~\cite{BZ}.

In the subsequent section, we use this degeneration to
prove that,
a homogeneous linear system of curves in $|dH|$ with $n$ general base points of
multiplicity $m$ is non-special if
all linear systems of curves in $|dH|$ with one multiple base point are
non-special.

Finally, in section~\ref{conj gamma=4}, we prove
conjecture~\cite[Conjecture~2.3~(i)]{DL}
for homogeneous linear systems of curves in $|dH|$ with $4^u9^w$ fat points
of multiplicity $m$,
if either $v \geq -1$ or $v \leq -1$ and $u >0$ or $2d \neq 1 \mod 3$.

\section{Preliminaries}

Let $S$ be a smooth projective generic $K3$ surface
and let $H$ be the generator of $\pic(S)$, then $H$ is ample, $H^2=2g-2 \geq 2$ and
$h^0(H)=g+1$;
moreover $H$ is very ample if $g \geq 3$
and if $g=2$ $H$ defines a double covering of $\p^2$ branched at an irreducible sextic
(see \cite[Proposition~3]{AM}).

Consider $Q_1,\cdots, Q_n$ points in general position on $S$,
for each one of
these points fix a multiplicity $m_1,\ldots, m_n$.
Define $\cl = \cl^{\gamma}(d,m_1,\ldots, m_r)$, with $\gamma = H^2 = 2g-2$,
as the linear system
of curves in $|dH|$ with multiplicity $m_i$ at $Q_i$ for all $i$.
Let $v$ denote $\vdim (\cl)$, then, using $ \dim |dH| =  \frac{\gamma d^2}{2} +1$,
we obtain that
\begin{equation}
v  = \frac{\gamma d^2}{2} +1 - \sum_{i=1}^n \frac{m_i(m_i + 1)}{2} \label{v = ..}
\end{equation}
Let $S'$ be the blowing-up of $S$ along the points $Q_1,\cdots, Q_n$,
$\pi : S' \map S$ the projection map and $E_i$ the exceptional divisor
corresponding to $Q_i$.
The linear system $\cl$ then corresponds to the system
$\pi^*(|dH|) \otimes \co_{S'}(-m_1 E_1 - \cdots -m_n E_n)$.
By abuse of notation, we will denote this linear system on the blowing-up
also by $\cl^{\gamma}(d,m_1,\ldots, m_r)$.

Similary, by $\cl^{\gamma}(d,m_1^{n_1},\ldots, m_y^{n_y})$,
we denote the linear system of curves in $|dH|$
with $n_i$ points (in general position) of multiplicity $m_i$ (for all $i = 1,\ldots,n$)
as well as the corresponding linear system on the blowing-up of $S$ along those
$n_1 + \cdots + n_y$ general points.

Let $Z=\sum_{i=1}^n m_i Q_i$ be the $0$-dimensional scheme defined by the multiple points
and consider the exact sequence of sheaves:
\[
\xymatrix@1{
0 \ar[r]  & \co_S(dH) \otimes {\mathcal I}_Z \ar[r]  & \co_S(dH) \ar[r]  & \co_Z \ar[r]  &  0
}
\]
where ${\mathcal I}_Z$ is the ideal sheaf of $Z$.
Taking cohomology we obtain
\begin{equation}
\label{equx}
v = h^0(\co_S(dH) \otimes {\mathcal I}(Z)) - h^1(\co_S(dH) \otimes {\mathcal I}(Z)) -1,
\end{equation}
because $h^1(\co_S(dH))=0$ (see e.g.~\cite{AM}).

Let $e$ denote the expected dimension of the linear system $\cl$,
i.e. $e  = \max\{ -1, v\}$
and let $l$ denote the dimension of $\cl$.
Then obviously we have that $ v \leq e \leq l$.

If $e < l$ (resp. $e = l$), we say the linear system $\cl$ is special (resp. non-special).
Note that equation~(\ref{equx}) shows that a non-empty linear system $\cl$
is non-special if and only if $h^1(\cl) = 0$.

\section{Degeneration of linear systems on $K3$ surfaces}\label{degen}

\subsection{The degeneration of a $K3$ surface}\label{degen surf}
~\\
Let $S$ be a $K3$ surface and $\Delta$ a complex disc around the origin.
Consider the product $V = S \times \Delta$
and its two projections
$q_1 : V \map S$ and $q_2 : V \map \Delta$.
Let $V_t$ denote $S \times \{t\}$.
Consider $b$ general points in $V_0$
and blow up $V$ along those $b$ points.
We then get a new threefold $X$ and the maps $\pi : X \map V$,
$p_1 = q_1 \circ \pi : X \map \Delta$
and $p_2 = q_2 \circ \pi : X \map S$;
i.e. we obtain the following commutative diagram
$$
\xymatrix{
      &   X  \ar[d]^{\pi} \ar@/_1pc/[ddl]_{p_1} \ar@/^1pc/[ddr]^{p_2} &     \\
      &   V  \ar[dl]_>>>>>>>>>{q_1} \ar[dr]^>>>>>>>>>{q_2}  &  \\
      \Delta  &  &  S
}
$$

Let $X_t$ be the fiber of $p_1$ over $t \in \Delta$.
If $t \neq 0$ then $X_t \cong V_t$ is our $K3$ surface $S$.
But, $X_0$ is the union of the proper transform $\st$ of $V_0$
and the $b$ exceptional divisors $\pp_i$.
Obviously each $\pp_i$ is isomorphic to $\p^2$
and $\st$ is the blowing-up of our $K3$ surface $S$ at
the $b$ general points
with projection map $\fb : \st \map S$.

Every $\pp_i$ intersects $\st$ transversally along a curve
$R_i$, which is a line in $\pp_i$ and an exceptional divisor on $\st$.
When we want to indicate that we consider $R_i$ in $\st$, resp. $\pp_i$,
we denote it by $E_i$, resp. $L_i$.

Note that the map $p_1$ gives a flat family of surfaces over $\Delta$,
so $X_0$ can be seen as a degeneration of $S$.

\subsection{The degeneration of a linear system on a $K3$ surface $S$}\label{degen sys}
~\\
Let $\cl$ be a line bundle on $S$,
and, for $k \in \z$, define the line bundle
$\co_{X}(\cl,k)$ on $X$ by
$$ \co_{X}(\cl,k) = p_2^* (\cl) \otimes \co_X(k\st). $$

The restriction of $ \co_{X}(\cl,k)$ to $X_t$, for $t \neq 0$, is then
isomorphic to $\cl$.

The restriction of $\co_{X}(\cl,k)$ to $X_0$,
which we denote by $\cx(\cl,k)$, is a flat limit of the line bundle $\cl$ on the general fiber $X_t$,
so $\cx(\cl,k)$ can be seen as a degeneration of the line bundle $\cl$.

On any $\pp_i$, we have that $\cx(\cl,k)|_{\pp_i} = \co_{X}(\cl,k)|_{\pp_i} \cong \co_{\p^2}(k)$,
since $\st$ intersects $\pp_i$ along a line.

On $\st$ we obtain that
$$\cx(\cl,k)|_{\st} = \co_{X}(\cl,k)|_{\st} = p_2^*(\cl)|_{\st} \otimes \co_X(k\st)|_{\st}.$$
But, since $\st \sim X_t - \sum_{i=1}^{l} \pp_i$ as divisors on $X$
and $\pp_i$ intersects $\st$ along $R_i$,
this means that
$$ \cx(\cl,k)|_{\st} \cong \fb^*(\cl) \otimes \co_{\st}(- \sum_{i=1}^{b} k E_i ). $$

Let $Q_1,\ldots,Q_n$ be general points on $S$, consider
a zero-dimensional subscheme $Z = m_1 Q_1 + \cdots + m_n Q_n$
and let $\cm$ be the sheaf $\cl \otimes \ci_Z$,
where $\ci_Z$ denotes the ideal sheaf that defines $Z$.
Choose positive integers $a_1, \ldots, a_b$ such that
$a_1 + \cdots a_b \leq n$.
Now, for all $i \in  \{1,\ldots,b\}$,
consider $a_i$ general points on $\pp_i$;
and take $n - \sum_{i=1}^b a_i$ general points on $\st$.
Denote those $n$ points (on the $\pp_i$ and $\st$) by $Q'_i$ (any order will do),
and let $Z'$ be the zero-dimensional subscheme of $X_0$ given by
$m_1 Q'_1 + \cdots + m_n Q'_n$.

Then we obtain that $\cx(\cl,k) \otimes \ci_{Z'}$
is a degeneration of $\cm$.

\subsection{Homogeneous linear systems on generic $K3$ surfaces}
~\\
Let $S$ be a generic $K3$ surface and consider
a homogeneous linear system $\cl = \cl^{\gamma}(d,m^n)$ on $S$.

Choose positive integers $b$ and $a$ such that $ab \leq n$,
let $X$ and $\cx(\cl,k)$ be as constructed before.
For all $1 \leq i \leq b$ and $1 \leq j \leq a$,
let $Q'_{i,j}$ be a general point on $\pp_i$
and for $1 \leq j \leq n -ab$
let $Q'_i$ be a general point on $\st$.
Now consider the zero-dimensional subscheme
$$ Z' = \sum_{i=1}^{n-ab} m Q'_i
    + \sum_{
        \begin{subarray}{c}
            i = 1,\ldots,b \\
            j = 1,\ldots,a
        \end{subarray}}
            m Q'_{i,j},
$$
On $X_0$ define $\cl_0  := \cx(\cl,k) \otimes \ci_{Z'}$,
$\cl_{i} := \cl_0|_{\pp_i}$ and $\cl_{\st} := \cl_0|_{\st}$.
Then
we see that
a divisor in the linear system
$\cl_0$ on $X_0$
consists of a divisor $D_{\st} \in \cl_{\st}$
and divisors $D_{i} \in \cl_{i}$
such that $D_{\st}|_{R_i} = D_i|_{R_i}$ for all $i$.

\begin{nots}
\begin{enumerate}
\item[] \svs
\item[] $l = \dim \cl$, $ l_0 = \dim \cl_0$, $ l_{\pp} = \dim \cl_i$, $ l_{\st} = \dim \cl_{\st}$ \svs
\item[] $ \crsi = \cl_{\st}|_{\cup_{i=1}^b E_i}$, $ r_{\st} = \dim \crsi$  \svs
\item[] $ \ccr_{i} = \cl_{\pp}|_{L_i}$, $ r_{\pp} = \dim \ccr_{i}$ \svs
\item[] $ \clh_{\st} = \cl_{\st} \otimes \co_{\st}(-\sum_{i=1}^b E_i)$, \svs
    $ \lh_{\st} = \dim \clh_{\st}$
\item[] $ \clh_{i} = \cl_{i} \otimes \co_{\pp_i}(- L_i)$, $ \lh_{\pp} = \dim \clh_{i}$
\end{enumerate}
\end{nots}

Obviously we have the following
\begin{align}
    l_{\st} &= r_{\st} + \lh_{\st} +1  \label{eq rs} \\
    l_{\pp} &= r_{\pp} + \lh_{\pp} +1. \label{eq rp}
\end{align}

Let $ \crpi $ denote the linear system on $\cup_{i=1}^b R_i$ for which $ \crpi|_{R_i} = \ccr_{i} $
and denote $r_{\cup L_i} = \dim \crpi$.
Then we have the following equality
\begin{equation} \label{eq dim}
l_0 = \dim (\crsi \cap \crpi) + b (\lh_{\pp} + 1) + \lh_{\st} + 1
\end{equation}
But, proceeding as in \cite[\S\thinspace 2]{BZ}, we obtain the transversality of
$\crsi$ and $\crpi$;
i.e.
\begin{align}
\dim (\crsi \cap \crpi) &= \max \{ -1 \, , \, r_{\st} + r_{\cup L_i} - b \, k \} \notag \\
    &= \max \{ -1 \, , \, r_{\st} + b \, r_{\pp} - b \, k \}.\label{eq inters}
\end{align}

For the virtual dimensions of the systems we introduce the following
\begin{nots}
\begin{enumerate}
\item[] \svs
\item[] $ v = \vdim \cl$ \svs
\item[] $ v_{\pp} = \vdim \cl_i$, $ v_{\st} = \vdim \cl_{\st}$ \svs
\item[] $ \vh_{\st} = \vdim \clh_{\st}$, $ \vh_{\pp} = \vdim \clh_{i}$
\end{enumerate}
\end{nots}

Since
$ \cl = \cl^{\gamma}(d,m^n) $,
$ \cl_{\st} \cong \cl^{\gamma}(d,k^b,m^{n-ab}) $,
$ \clh_{\st} \cong \cl^{\gamma}(d,(k+1)^b,m^{n-ab}) $,
$ \cl_{i} \cong \cl (k, m^a)$
and $\clh_i \cong \cl (k-1 , m^a)$;
a simple calculation shows that
\begin{align}
    v   &=  v_{\st}  +  b \, \vh_{\pp} + b
        \, =  v_{\st} + b \, (v_{\pp} - k )   \notag \\
        &=  \vh_{\st}  + b \, v_{\pp} + b \label{eq vdim}
        \, =  \vh_{\st} + b \, (\vh_{\pp} + k +2)
\end{align}

\subsection{Remark}\label{double degen}
~\\
Let $S$ be a $K3$ surface, $\cl$ a line bundle on $S$
and let $X$ and $\cx(\cl,k)$ be as before.

Now do the construction of \S\thinspace\ref{degen surf} and \S\thinspace\ref{degen sys}
using $X_0$ in stead of $S$ and
$\cx(\cl,k)$ (or $\cx(\cl,k) \otimes \ci_{Z'}$) in stead of $\cl$ (or $\cl \otimes \ci_Z$);
i.e. consider $W = X_0 \times \Delta$, blow up $W$ along $b'$ general points on
$\st \times \{0\} \subset W_0 = X_0 \times \{0\}$ and obtain the following commutative diagram
$$
\xymatrix{
      &   Y  \ar[d]^{\pi'} \ar@/_1pc/[ddl]_{p'_1} \ar@/^1pc/[ddr]^{p'_2} &     \\
      &   W  \ar[dl]_>>>>>>>>>{q'_1} \ar[dr]^>>>>>>>>>{q'_2}  &  \\
      \Delta  &  &  X_0
}
$$

So we obtain a degeneration of $X_0$, and
a degeneration of
$\cx(\cl,k)$ (or $\cx(\cl,k) \otimes \ci_{Z'}$).

We call this a {\em double degeneration}
of $S$ and $\cl$ (or $\cl \otimes \ci_Z$),
and, continuing in the same way, we can obtain an $\eta$-uple
degeneration of $S$ and $\cl$ (or $\cl \otimes \ci_Z$), for any $\eta \geq 2$.

\section{Applying the degeneration to homogeneous linear systems on generic $K3$ surfaces}
\label{degen generic}

From now on we assume that we work on a generic $K3$ surface $S$,
with $H^2 = \gamma$, where $H$ is the generator of $\pic S$.

The main result of this section is
the following

\begin{theor}\label{4^u 9^w}
If $\cl^{\gamma}(d,\mu)$ is non-special for all $\mu$,
then $\cl = \cl^{\gamma}(d,m^{n})$ with $n = 4^u 9^w$ is non-special
for all positive integers $m, u$ and $w$.
\end{theor}

To make the proof of this theorem more transparent we first fix some notation
and state a few auxiliary results.

Let $c \in \{4,9\}$ such that $c | n$ and consider the degeneration $(X_0 , \cl_0)$
of $(S , \cl)$ obtained by using the construction explained in \S\thinspace\ref{degen} with
$b = n/c$ and $k \in \n$.
Note that, since $\cl_0$ is a degeneration of $\cl$,
$l_0 \geq l \geq e \geq v$.

On $\p^2$, let $\cl(\delta,\mu^{\nu})$ denote the linear system of
plane curves of degree $\delta$ having multiplicity $\mu$
at $\nu$ points (in general position).

The following is a well know result, and can easily be checked using for instance
the results of~\cite{BH}.

\begin{lemma}\label{4 and 9}
On $\p^2$ the linear system $\cl(\delta,\mu^c)$ with $c \in \{4,9\}$
is non-special for all $\delta$ and $\mu$. \vrk
\end{lemma}

\begin{rem}
Lemma~\ref{4 and 9} implies in particular that $\cl_i$ and $\clh_i$
are non-special linear systems on $\pp_i$.
\end{rem}

\begin{claim}\label{max >= -1}
If $v \geq -1$,
$v_{\st} \geq -1$, $v_{\pp} \geq -1$
and $\cl_{\st}$ and $\clh_{\st}$ are non-special systems;
then $\dim (\crsi \cap \crpi) = r_{\st} + b \, r_{\pp} - b \, k$.
\end{claim}

\begin{proof}
Note that the conditions of the theorem immediately imply that
$l_{\st} = v_{\st}$ and $l_{\pp} = v_{\pp}$.

Since $\dim (\crsi \cap \crpi) = \max\{ -1 , r_{\st} + b \, r_{\pp} - b \, k \}$,
it suffices to show that
$r_{\st} + b \, r_{\pp} - b \, k \geq -1$.

Because of (\ref{eq rs}) and (\ref{eq rp}),
we have that
$$
r_{\st} + b \, r_{\pp} - b \, k
        = l_{\st} - \lh_{\st} -1 + b \,( l_{\pp} - \lh_{\pp} - 1) - b \, k.
$$

If $\vh_{\st} \geq -1$ and $\vh_{\pp} \geq -1$,
then $\lh_{\st} = \vh_{\st}$ and $\lh_{\pp} = \vh_{\pp}$.
So
\begin{align*}
    r_{\st} + b \, r_{\pp} - b \, k
        &= v_{\st} - \vh_{\st} -1 + b \,( v_{\pp} - \vh_{\pp} - 1) - b \, k \\
        &= b \, (k + 1) - 1 \geq -1.
\end{align*}

If $\vh_{\st} \leq -2$ and $\vh_{\pp} \geq -1$,
then $\lh_{\st} = -1$ and $\lh_{\pp} = \vh_{\pp}$.
So
\begin{align*}
    r_{\st} + b \, r_{\pp} - b \, k
        &= v_{\st} + b \,( v_{\pp} - \vh_{\pp} - 1) - b \, k \\
        &= v_{\st} \geq -1.
\end{align*}

If $\vh_{\st} \geq -1$ and $\vh_{\pp} \leq -2$,
then $\lh_{\st} = \vh_{\st}$ and $\lh_{\pp} = -1$.
So
\begin{align*}
    r_{\st} + b \, r_{\pp} - b \, k
        &= v_{\st} - \vh_{\st} -1 + b \, v_{\pp}  - b \, k \\
        &= b \, (1 + v_{\pp}) - 1  \geq -1.
\end{align*}

If $\vh_{\st} \leq -2$ and $\vh_{\pp} \leq -2$,
then $\lh_{\st} = -1$ and $\lh_{\pp} = -1$.
So
\begin{align*}
    r_{\st} + b \, r_{\pp} - b \, k
        &= v_{\st} + b \, (v_{\pp} - k) \\
        &\!\!\overset{\mbox{{\rm \tiny{(\ref{eq vdim})}}}}{=} v  \geq -1.
\end{align*}
\end{proof}

\begin{lemma}\label{v >=- 1, ind step}
If $v \geq -1$,
$v_{\st} \geq -1$, $v_{\pp} \geq -1$
and $\cl_{\st}$ and $\clh_{\st}$ are non-special systems;
then $\cl$ is non-special.
\end{lemma}

\begin{proof}
Because $\cl_0$ is a degeneration of $\cl$,
we know that $v \leq l \leq l_0$;
so it suffices to prove that $l_0 = v$.

Using (\ref{eq dim}), (\ref{eq inters}) and claim~\ref{max >= -1}
we obtain that
\begin{equation*}
l_0 = r_{\st} + b \, r_{\pp} - b \, k + b \, (\lh_{\pp}+1) + \lh_{\st} +1.
\end{equation*}

So, using (\ref{eq rs}) and (\ref{eq rp}),
we see that
\begin{equation*}
l_0 = l_{\st} + b \, (l_{\pp} - k)
    = v_{\st} + b \, (v_{\pp} - k)
    \overset{\mbox{{\rm \tiny{(\ref{eq vdim})}}}}{=} v.
\end{equation*}
\end{proof}

\begin{lemma}\label{v >=- 1, k exists}
Let $\cl = \cl^{\gamma}(d,m^{n})$ with $n = 4^u 9^w$,
$d, m , u, w \in \n$, $ u+w > 0$ and $v \geq -1$.
Take $c$, $b$, $X_0$ and $\cl_0$ as before.
Then $\exists \, k \in \n$ such that
            $v_{\st} \geq -1$ and $v_{\pp} \geq -1$.
\end{lemma}

\begin{proof}
Because
$v_{\st} = \frac{\gamma}{2}d^2 +1 - b \frac{k(k+1)}{2}$
and
$v_{\pp} = \frac{k(k+3)}{2} - c \frac{m(m+1)}{2}$,
the inequalities
$v_{\st} \geq -1$ and $v_{\pp} \geq -1$
are equivalent to
\begin{align}
 k^2 + k - \alpha \leq 0 , \ \ \ &\mbox{ with } \alpha = \frac{1}{b}(\gamma d^2 +4) \ \mbox{and} \\
 k^2 + 3k - \beta \geq 0 , \ \ \ &\mbox{ with } \beta  = c m(m+1) - 2.
\end{align}
But this is the same as
\begin{align}
k \in  &\left[ \frac{-1-\sqrt{1+4 \alpha}}{2}, \frac{-1+\sqrt{1+4 \alpha}}{2} \right] \ \mbox{and}  \label{eq1}\\
k \in  & \, \left] -\infty,\frac{-3-\sqrt{9+4 \beta}}{2}\right] \cup
        \left[ \frac{-3+\sqrt{9+4 \beta}}{2}, \infty \right[.\label{eq2}
\end{align}
So proving the statement is equivalent to proving that
there exists a positive integer $k$ such that both
(\ref{eq1}) and (\ref{eq2}) are satisfied, i.e.
it is enough to show that
\begin{equation*}
    \frac{-3+\sqrt{9+4 \beta}}{2} + 1 \leq \frac{-1+\sqrt{1+4 \alpha}}{2}.
\end{equation*}
A simple calculation shows that the previous inequality is equivalent to
$ \alpha \geq \beta +2$,
which is in turn equivalent to $v \geq -1$.
\end{proof}

\begin{lemma}\label{v <=- 1, ind step}
If $v \leq -1$,
$\vh_{\st} \leq -1$, $\vh_{\pp} \leq -1$
and $\cl_{\st}$ and $\clh_{\st}$ are non-special systems;
then $\cl$ is non-special and thus empty.
\end{lemma}

\begin{proof}
Since $v_{\st} = \vh_{\st} + b(k+1)$ and $v_{\pp} = \vh_{\pp} + k +1$,
we see that $v_{\st} \leq l(k+1)-1$ and $v_{\pp} \leq k $.

If $v_{\st} \leq -1$ and $v_{\pp} \leq -1$
then $r_{\st} = r_{\pp} = -1$.
If $v_{\st} \leq -1$ and $v_{\pp} \geq -1$,
then $r_{\st} = -1$ and $r_{\pp} = v_{\pp} \leq k$.
If $v_{\st} \geq -1$ and $v_{\pp} \leq -1$,
then $r_{\st} = v_{\st} \leq b (k+1)-1$ and $r_{\pp} = -1$.
If $v_{\st} \geq -1$ and $v_{\pp} \leq -1$,
then $r_{\st} = v_{\st}$ and $r_{\pp} =  v_{\pp}$.
So in any case we obtain that
$r_{\st} + b \, r_{\pp} - b \, k \leq -1$,
i.e. $\dim (\crsi \cap \crpi) = -1$.

Using (\ref{eq dim}) and $\lh_{\st} = \lh_{\pp} = -1$,
we see that
$l_0 = \lh_{\st} + k (\lh_{\pp} +1) = -1$.
\end{proof}

\begin{lemma}\label{v <=- 1, k exists}
Let $\cl = \cl^{\gamma}(d,m^{n})$ with $n = 4^u 9^w$,
$d, m , u, w \in \n$, $ u+w > 0$ and $v \leq -1$.
Take $c$, $b$, $X_0$ and $\cl_0$ as before.
Then $\exists \, k \in \n$ such that
            $\vh_{\st} \leq -1$ and $\vh_{\pp} \leq -1$.
\end{lemma}

\begin{proof}
Because
$\vh_{\st} = \frac{\gamma}{2}d^2 +1 - b \frac{(k+1)(k+2)}{2}$
and
$\vh_{\pp} = \frac{(k-1)(k+2)}{2} - c \frac{m(m+1)}{2}$,
the inequalities
$\vh_{\st} \leq -1$ and $\vh_{\pp} \leq -1$
are equivalent to
\begin{align*}
 k^2 + 3k - \alpha \geq 0 , \ \ \ &\mbox{ with } \alpha = \frac{1}{b}(\gamma d^2 +4) - 2 \ \mbox{and} \\
 k^2 + k - \beta \leq 0 , \ \ \ &\mbox{ with } \beta  = c m(m+1).
\end{align*}
Proceeding as in the proof of lemma~\ref{v >=- 1, k exists},
we obtain that it is sufficient to prove that
\begin{equation*}
    \frac{-1+\sqrt{1+4 \beta}}{2} \geq \frac{-3+\sqrt{9+4 \alpha}}{2} +1,
\end{equation*}
which is equivalent to $ \beta \geq \alpha +2$.
And this last inequality is equivalent to $v \geq -1$.
\end{proof}

\begin{proof}[Proof of theorem~\ref{4^u 9^w}]
~\\
Using induction, the result follows immediately from
lemma's~\ref{v >=- 1, ind step} and~\ref{v >=- 1, k exists} if $v \geq -1$,
and from lemma's~\ref{v <=- 1, ind step} and~\ref{v <=- 1, k exists} if $v \leq -1$.
\end{proof}

\subsection{Using a higher order degeneration}
~\\
Using the $\eta$-uple degeneration explained in Remark\ref{double degen},
we can, proceeding as above, prove that
the non-speciality of $\cl^{\gamma}(d,\mu_1,\ldots, \mu_{\eta})$
for all $\mu_1,\ldots,\mu_{\eta}$ implies the non-speciality of
$\cl = \cl^{\gamma}(d,m_1^{n_1},\ldots, m_{\eta}^{n_\eta})$
with $n_i=4^{u_i}9^{w_i}$
for all $m_1,\ldots,m_{\eta}, u_1,\ldots, u_{\eta}, w_1,\ldots,w_{\eta}$.

As we will see in section~\ref{1 mult},
proving the non-speciality of systems $\cl^{\gamma}(d,\mu)$
is already rather complex.
In fact, in section~\ref{1 mult} we only prove
when such a system is non-special for
general generic $K3$ surfaces with $\gamma = 4$.

\subsection{Remark}
~\\
As proved in~\cite[Theorem~6.1]{CM2}, the Segre conjecture on planar linear systems through fat points
implies that the only special homogeneous linear systems $\cl(k, m^n)$ are the ones
with $n \in \{2,3,5,6,7,8\}$.
So if the Segre conjecture is true,
we can do the degeneration using, not only $a=4$ or $9$, but
$a \in \{4\} \cup \z_{\geq 9}$.
In this way we can then
prove theorem~\ref{4^u 9^w},
for any $n$ which can be written as a product of powers of numbers in $\{4\} \cup \z_{\geq 9}$,
and thus in particular for any $n \geq 9$.

\section{The non-speciality of linear systems $\cl^{\gamma}(d,m^{n})$ with $n = 4^u 9^w$
    on general generic $K3$ surfaces with $\gamma = 4$
    }\label{conj gamma=4}

The main result of this section is the following

\begin{theor}\label{gamma=4}
Let $S$ be a general generic $K3$ surface
with $\gamma = 4$
and let
$\cl = \cl^{4}(d,m^n)$, with
$n = 4^u 9^w$, $u,w \in \z_{\geq 0}$ and $d \in \z_{>0}$.
\svs \\
\ind (1) If $v(\cl) \geq -1$ then the linear system $\cl$ is non-special.
\svs \\
\ind (2) If $v(\cl) \leq -1$ and either $u > 0$ or $ 2d \neq 1 \mod 3$
        then $\cl$ is non-special unless $n = 1$, $m = 2d$ and $d \geq 2$.
In the latter case $\dim (\cl^{4}(d,2d)) = 0 > \edim (\cl^{4}(d,2d)) = -1$
and $\cl^{4}(d,2d) = d C$, with $C$ the unique element of $\cl^4(1,2)$.
\end{theor}

\begin{rem}
Note that the previous theorem implies that \cite[Conjecture~2.3~(i)]{DL}
is true
for homogeneous linear systems with $n = 4^u 9^w$ fat points on $S$
in the following cases
\svs \\
\ind (i) $v \geq -1$, or
\svs \\
\ind (ii) $v \leq -1$ and $u > 0$, or
\svs \\
\ind (iii) $v \leq -1$ and $2d \neq 1 \mod 3$.
\end{rem}

We will prove theorem~\ref{gamma=4} by using the degeneration as introduced
in \S\thinspace\ref{degen generic}
and the following

\begin{prop}\label{1 mult}
Let $S$ be a general generic $K3$ surface
with $\gamma = 4$.
Then $\cl = \cl^{4}(d,\mu)$, with $d, \mu \in \z_{>0}$, is non-special
unless $\mu = 2d$ and $d \geq 2$.
In the latter case $\dim (\cl^{4}(d,2d)) = 0 > \vdim (\cl^{4}(d,2d)) = -1$
and $\cl^{4}(d,2d) = d C$, with $C$ the unique element of $\cl^4(1,2)$.
\end{prop}

\begin{proof}
The map $\phi : S \map \p^3$ corresponding to $|H|$ (with $H$ the generator of $\pic S$)
is an embedding, so we may look at $S$ as a quartic surface in $\p^3$.
The unique element $C$ of $\cl^4(1,2)$ is then the divisor on $S \subset \p^3$
defined by the tangent plane to $S$ at this general point $P$.
So $C$ is an irreducible plane curve of degree 4 having a node at $P$ (and no other singularities).

Blow up $S$ along the point $P$, and, by abuse of notation, let $\cl^4(d,\mu)$ also denote
the complete linear system on this blowing-up corresponding to $\cl^4(d,\mu)$ on $S$.
Let $\ct$ denote the strict transform of $C$ on this blowing-up,
then $g(\ct) = 2$.

Assume $\mu = 2d$ ($d \geq 2$) and consider the following sequence
$$
\xymatrix@1{
  0  \ar[r]  &  \cl^4(d-1,2d-2)  \ar[r]  &   \cl^4(d,2d)  \ar[r]  &  \cl^4(d,2d) \otimes \co_{\ct}
     \ar[r]  &  0.
}$$
Because $\deg ( \cl^4(d,2d) \otimes \co_{\ct}) = 0 = 2g(\ct) - 4$,
we know that $h^0(\cl^4(d,2d) \otimes \co_{\ct}) = 0$
unless $\cl^4(d,2d) \otimes \co_{\ct} = \co_{\ct}$.
But $\cl^4(d,2d) \otimes \co_{\ct} = | d g^2_4 - 2d (P_1 + P_2)|$,
where $P_1$ and $P_2$ are the intersection points of $\ct$ with the exceptional curve.
So this would mean that
$| d g^2_4 - 2d (P_1 + P_2)| = \co_{\ct}$.
But, since $P$ is a general point on the general quaritc $S$,
we know that $g^2_4 = |K_{\ct} + P_1+P_2|$;
so we would obtain that 
$| d K_{\ct} - d (P_1 + P_2)| = \co_{\ct}$,
or thus that $| d K_{\ct} | = |d (P_1 + P_2)|$,
which is not true (because $P$ is a general element of $S$).
So $h^0(\cl^4(d,2d)) = h^0(\cl^4(d-1,2(d-1)))$ for all $d \geq 2$.
Using this a number of times we thus obtain that
$h^0(\cl^4(d,2d)) = h^0(\cl^4(1,2)) = 1$,
so $d C$ is the only divisor in $\cl^4(d,2d)$.

It now follows immediately that $\dim \cl = -1$ if $\mu \geq 2d+1$,
since the only divisor in $\cl^4(d,2d)$ has multiplicity exactly $2d$ in $P$.

Now consider the case $\mu \leq 2d-1$.
If we can prove that $\cl$ is non-special for $\mu = 2d-1$, then
the non-speciality follows for
all $\mu \leq 2d-1$ (since $\vdim \cl^4(d,2d-1) = d+1$).
For $d=1$, there is nothing to prove, since $P$ is a general point of $S$.
So we may assume that $d \geq 2$ and that the
$\cl^4(d',2d'-1)$ is non-special for $d' \leq d-1$.
Now consider the following sequence
$$
\xymatrix@1{
  0  \ar[r]  &  \cl^4(d-1,2d-3)  \ar[r]  &   \cl^4(d,2d-1)  \ar[r]  &  \cl^4(d,2d-1) \otimes \co_{\ct}
     \ar[r]  &  0.
}$$
Because $\deg ( \cl^4(d,2d-1) \otimes \co_{\ct}) = 2 = 2g(\ct) - 2$,
we know that $h^1(\cl^4(d,2d-1) \otimes \co_{\ct}) = 0$
unless $\cl^4(d,2d-1) \otimes \co_{\ct} = K_{\ct}$ (the canonical class on $\ct$).
But $\cl^4(d,2d) \otimes \co_{\ct} = | d g^2_4 - (2d-1) (P_1 + P_2)|$
and $g^2_4 = |K_{\ct} + P_1+P_2|$
(where $P_1$ and $P_2$ are, as before, the intersection points of $\ct$ with the exceptional curve).
Then this would mean that
$| (d-1) K_{\ct}| = |(d - 1)(P_1 + P_2)| $,
which is not true (because $P$ is a general element of $S$).
So $h^1( \cl^4(d,2d-1) \otimes \co_{\ct})=0$ and, by hypotheses, $h^1(\cl^4(d-1,2d-3))=0$,
thus also $h^1( \cl^4(d,2d-1) ) =0$.
\end{proof}

\begin{proof}[Proof of theorem~\ref{gamma=4}]\svs
~\svs \\
(1)
In case $n=1$ we are done because of proposition~\ref{1 mult}.
So we assume $n > 1$, or equivalently $u+w \geq 1$.
As in the proof of theorem~\ref{4^u 9^w}, we can use
lemma's~\ref{v >=- 1, ind step} and \ref{v >=- 1, k exists}, but we have to
end the induction with a class $\cl_{\st}=\cl^4(d,\mu)$ such that $\mu \neq 2d$ and $\mu \neq 2d-1$
(in order to apply lemma~\ref{v >=- 1, ind step}
we need to have $\cl_{\st}$ and $\clh_{\st}$ non-special).
So, it is enough to prove that, in case $ u+w =1$, we can take
$k$ such that $v_{\st} \geq -1$, $v_{\pp} \geq -1$ and
$k \notin \{2d,2d-1\}$.
Since now $v_{\st} \geq -1$, is equivalent to $ k \leq 2d-1$,
it is enough to prove that $v_{\pp} \geq -1$ for $k=2d-2$
(which then implies that you can always take $k=2d-2$ in the last step).
Now, $v_{\pp} = \frac{k(k+3)}{2} - c \frac{m(m+1)}{2}$,
so, for $k=2d-2$, we obtain
$v_{\pp} = 2d^2 - d - 1 - c \frac{m(m+1)}{2}$.
Since $ u+w =1$, we have that $c = n$, so $c \frac{m(m+1)}{2} = 2d^2 +1 - v$,
and we see that $v_{\pp} \geq -1$ for $k=2d-2$
if and only if $ v \geq d+1$.
Now, we only know that $v \geq -1$, and $c = n \in \{4,9\}$,
but a simple calculation shows that
in case $n=4$:
$$v \geq -1 \iff m \leq d-1  \iff v \geq 1+2d. $$
In case $n=9$, we may assume that $d \geq 4$,
since for $d \leq 3$, we obtain that $v \geq -1$ is only possible if $m \leq 1$.
And, a simple calculation then shows that
$$ v \geq -1  \iff m \leq \frac{2d-2}{3}  \iff v \geq 2+d$$
(note that, if $m = \frac{2d-1}{3}$, we obtain $v = 2-d < -1$, since $d \geq 4$).
\svs
\\
(2)
Again, if $n=1$ we are done because of proposition~\ref{1 mult}.
So we assume $n > 1$, or equivalently $u+w \geq 1$.
As in the proof of theorem~\ref{4^u 9^w}, we can use
lemma's~\ref{v <=- 1, ind step} and \ref{v <=- 1, k exists},
up to the last step.
For this last step, we have $c = n$ and $b=1$, and we will take $k=2d$.
In this case, we can no longer use lemma~\ref{v <=- 1, ind step} since
$\cl_{\st}$ is special.
If we can prove however that
$\dim (\crsi \cap \crpi) = \lh_{\st} = \lh_{\pp} = -1$,
then we still obtain $l_0 =-1$, and thus also $l = -1$.
Since $\dim (\crsi \cap \crpi) = \max \{-1, r_{\st} +  (r_{\pp} - 2d)\}$,
$r_{\st} = 0$ and $r_{\pp} = l_{\pp} $
it is enough to prove that, for $k=2d$,
$ v_{\pp} \leq 2d-1$ and $\vh_{\pp} \leq -1$.
These inequalities are equivalent to respectively
$v \leq -d$ and $v \leq 1 -d$,
so it is sufficient to have $v \leq -d$.

In case $u >0$, we can make sure that in the last step $c = 4$,
and we obtain
$$ v \leq -1  \iff  m \geq d  \iff  v \leq 1-2d. $$
In case $u =0$, we have $c = 9$ in the last step, and
using the fact that $2d \neq 1 \mod 3$, a simple calculation shows that
$$ v \leq -1  \iff  m \geq \frac{2d}{3}  \iff  v \leq 1-3d. $$
\end{proof}

\section*{Acknowledgement}

Both authors would like to thank Marc Coppens for helpfull discussions on the subject.



\begin{thebibliography}{May72}

\bibitem[BZ03]{BZ}
Anita Buckley and Marina Zompatori.
\newblock Linear systems of plane curves with a composite number of base points
  of equal multiplicity.
\newblock {\em Trans. Amer. Math. Soc.}, 355(2):539--549 (electronic), 2003.

\bibitem[CM98]{CM1}
Ciro Ciliberto and Rick Miranda.
\newblock Degenerations of planar linear systems.
\newblock {\em J. Reine Angew. Math.}, 501:191--220, 1998.

\bibitem[CM01]{CM2}
Ciro Ciliberto and Rick Miranda.
\newblock The {S}egre and {H}arbourne-{H}irschowitz conjectures.
\newblock In {\em Applications of algebraic geometry to coding theory, physics
  and computation (Eilat, 2001)}, volume~36 of {\em NATO Sci. Ser. II Math.
  Phys. Chem.}, pages 37--51. Kluwer Acad. Publ., Dordrecht, 2001.

\bibitem[DL03]{DL}
Cindy De Volder and Antonio Laface.
\newblock Linear systems on generic $K3$ surfaces.
\newblock {\em Preprint, math.AG/0309073}, 2003.

\bibitem[Gim89]{AG}
Alessandro Gimigliano.
\newblock Regularity of linear systems of plane curves.
\newblock {\em J. Algebra}, 124(2):447--460, 1989.

\bibitem[Har85]{BH}
Brian Harbourne.
\newblock Complete linear systems on rational surfaces.
\newblock {\em Trans. Amer. Math. Soc.}, 289(1):213--226, 1985.

\bibitem[May72]{AM}
Alan~L. Mayer.
\newblock Families of {$K3$} surfaces.
\newblock {\em Nagoya Math. J.}, 48:1--17, 1972.

\end{thebibliography}

\end{document}